\documentclass[11pt]{article}
\usepackage[OT2,T1]{fontenc}
\usepackage{amsmath}\usepackage{amssymb}\usepackage{amscd}
\newcommand\cyr{%
\renewcommand\rmdefault{wncyr}%
\renewcommand\sfdefault{wncyss}%
\renewcommand\encodingdefault{OT2}%
\normalfont\selectfont}
\DeclareTextFontCommand{\textcyr}{\cyr}

\DeclareSymbolFont{cyrletters}{OT2}{wncyr}{m}{n}
\DeclareSymbolFont{cyrlettersb}{OT2}{wncyr}{b}{n}
\DeclareMathSymbol{\Sha}{\mathalpha}{cyrletters}{"58}
\DeclareMathSymbol{\sha}{\mathalpha}{cyrlettersb}{"78}

\def\be{\begin{equation}}\def\ee{\end{equation}}\def\ba{\begin{eqnarray}}\def\ea{\end{eqnarray}}
\def\lb{\label}
\def\s{\sigma}
\def\lb{\label}\def\be{\begin{equation}}\def\ee{\end{equation}}
\newcommand{\binc}[2]{\Bigl(\!\begin{array}{c}{#1}\\ {#2}\end{array}\!\Bigr)}

\begin{document}
\begin{center}
\hfill {\small{CPT-P094-2008}}

\bigskip
{\Large\bf Braids, Shuffles and Symmetrizers}

\vspace{.6cm}

{\large\bf A. P. Isaev}
 
\vskip .2cm

Bogoliubov Laboratory of Theoretical Physics, JINR\\
141980 Dubna, Moscow Region, Russia

\vskip .4cm

{\large {\bf O. V. Ogievetsky\footnote{On leave of absence from P. N.
Lebedev Physical Institute, Leninsky Pr. 53, 117924 Moscow, Russia}}}

\vskip 0.2 cm

Center of Theoretical Physics\footnote{Unit\'e Mixte de Recherche (UMR 6207) du CNRS et 
des Universit\'es Aix--Marseille I, Aix--Marseille II et du Sud Toulon -- Var; laboratoire 
affili\'e \`a la FRUMAM (FR 2291)}, Luminy \\
13288 Marseille, France
\end{center}

\vskip .8cm
\hfill\textcyr{Pamyati Aleshi Zamolodchikova\hspace{1cm}}

\vspace{1.0cm}

\centerline{\footnotesize\bf Abstract} 
{\footnotesize\noindent 
\begin{center}\begin{minipage}[c]{14cm}
Multiplicative analogues of the shuffle elements of the braid group rings are introduced; 
in local representations they give rise to certain graded associative algebras (b-shuffle 
algebras). For the Hecke and BMW algebras, the (anti)-symmetrizers have simple expressions 
in terms of the multiplicative shuffles. The (anti)-symmetrizers can be expressed in terms 
of the highest multiplicative 1-shuffles (for the Hecke and BMW algebras) and in terms of 
the highest additive 1-shuffles (for the Hecke algebras). The spectra and multiplicities of 
eigenvalues of the operators of the multiplication by the multiplicative and additive 
1-shuffles are examined. 
\end{minipage}\end{center}}

\vspace{.1 cm}
\tableofcontents

\newpage
\section{Braid shuffles}

In this section we collect some necessary information on shuffle elements in the braid group rings.

\vskip .2cm
In the Artin presentation, the braid group $B_{M+1}$ is given by generators 
$\sigma_i$, $1\leq i\leq M$, and relations 
\begin{eqnarray}\lb{ryb1}\sigma_{i}\sigma_{j}\sigma_{i}=\sigma_{j}\sigma_{i}\sigma_{j}
&{\mathrm{if}}& |i-j|=1\ ,\\[.5em]
\sigma_{i}\sigma_{j}=\sigma_{j}\sigma_{i}&{\mathrm{if}}&|i-j|>1\ .\label{braidg}\end{eqnarray}
The inductive limit $B_\infty =\displaystyle{\lim_{\longrightarrow}}\, B_M$ is defined by 
inclusions $B_M\rightarrow B_{M+1}$, 
$B_M\ni\sigma_i\mapsto \sigma_i\in B_{M+1}$, $i=1,\dots ,M-1$. 

\vskip .2cm
We denote $w^{\uparrow\ell}$, as in \cite{OP}, the image of an element $w\in B_\infty$ under
the endomorphism of $B_\infty$, sending $\sigma_i$ to $\sigma_{i+\ell}$, $i=1,2,\dots$
(we keep the same notation for the Hecke and BMW quotients of the braid group rings). 

\vskip .2cm
Braid shuffle elements $\Sha_{m,n}$ ($m,n\in {\mathbb Z}_{\geq 0}$) are analogues of the 
binomial coefficients. The shuffle elements belong to the group ring of $B_{m+n}$ 
(and thereby of $B_\infty$); they can be defined inductively by any of the recurrence relations 
(braid analogues of the Pascal rule)
\begin{eqnarray}\lb{shufrec}\Sha_{m,n}^{\phantom{\uparrow}}&=&\Sha_{m-1,n}^{\phantom{\uparrow}}+
\Sha_{m,n-1}^{\phantom{\uparrow}}\sigma_{m+n-1}\cdots\sigma_{n}\ ,\\[.5em]
\lb{shufrec2}\Sha_{m,n}^{\phantom{\uparrow}}&=&
\Sha_{m,n-1}^{\uparrow 1}+\Sha_{m-1,n}^{\uparrow 1}\sigma_{1}\cdots\sigma_{n}\ ,\end{eqnarray} 
together with the boundary conditions $\Sha_{0,n}=1$ and $\Sha_{n,0}=1$ for any non-negative
integer $n$.

\vskip .2cm
Let $\Sigma_n$ be the lift \cite{Ma} of the symmetrizer $\sum_{g\in S_n}g$ from the symmetric group ring ${\mathbb Z}{\mathbb{S}}_n$ to ${\mathbb Z}B_n$. The element $\Sigma_n$ is the 
braid analogue of $n!$; it satisfies
\be\lb{symsh}\Sigma_{m+n}^{\phantom{\uparrow}}=\Sha_{n,m}^{\phantom{\uparrow}}\,
\Sigma_m^{\phantom{\uparrow}}\,\Sigma_n^{\uparrow m}\ .\ee
Using the automorphism $\mathfrak a$ and the anti-automorphism $\mathfrak b$, 
${\mathfrak{b}}(xy)={\mathfrak{b}}(y){\mathfrak{b}}(x)$,
of the braid group 
$B_{n+1}$, defined on the generators by
\be\lb{aau}\mathfrak a:\s_i\mapsto\s_{n+1-i}\quad ,\quad\mathfrak b:\s_i\mapsto\s_i\ ,\ee
and their composition, one obtains three more decompositions of $\Sigma$.

\vskip .2cm
Higher shuffles (braid analogues of the trinomial {\em etc.} coefficients) appear in the
further decompositions of the elements $\Sigma_n$,
\be\lb{hsh}\Sigma_{m+n+k}^{\phantom{\uparrow}}=\left\{
\begin{array}{l}\Sha_{n+k,m}^{\phantom{\uparrow}}\,
\Sigma_m^{\phantom{\uparrow}}\,\Sigma_{n+k}^{\uparrow m}
=\Sha_{n+k,m}^{\phantom{\uparrow}}\,\Sha_{k,n}^{\uparrow m}\,\Sigma_m^{\phantom{\uparrow}}
\Sigma_n^{\uparrow m}\,\Sigma_{k}^{\uparrow m+n}\ ,\\[1.2em]
\Sha_{k,m+n}^{\phantom{\uparrow}}\,\Sigma_{m+n}^{\phantom{\uparrow}}\,\Sigma_{k}^{\uparrow m+n}
=\Sha_{k,m+n}^{\phantom{\uparrow}}\,\Sha_{n,m}^{\phantom{\uparrow}}\,\Sigma_m^{\phantom{\uparrow}}
\Sigma_n^{\uparrow m}\,\Sigma_{k}^{\uparrow m+n}\ .\end{array}\right.\ee
Due to the existence of a (one-sided) order on the braid groups \cite{D}, the braid group rings 
${\mathbb Z}B_n$ have no zero divisors. Equating the two expressions for $\Sigma_{m+n+k}$ 
in (\ref{hsh}) and simplifying, one finds
\be\lb{ctsh}\Sha_{n+k,m}^{\phantom{\uparrow}}\,\Sha_{k,n}^{\uparrow m}=
\Sha_{k,m+n}^{\phantom{\uparrow}}\,\Sha_{n,m}^{\phantom{\uparrow}}\ .\ee
A direct verification of (\ref{ctsh}) is a good exercise. Any of the expressions in (\ref{ctsh})
is the braid trinomial coefficient $\Sha_{k,n,m}$. The element $\Sigma_n$ is the shuffle
$\Sha_{1,1,\dots ,1}$.

\vskip .2cm
We shall later use the following identity
\be\lb{lue}\Sha_{1,n-1}^{\phantom{\uparrow}}\Sha_{1,n-2}^{\phantom{\uparrow}}\dots
\Sha_{1,n-k}^{\phantom{\uparrow}}=
\Sha_{k,n-k}^{\phantom{\uparrow}}\Sigma_{k}^{\uparrow n-k}\ ,\ee
which is verified by induction. For $k=1$ there is nothing to prove. The induction step
uses (\ref{ctsh}) and then (\ref{symsh}):
\be\lb{neid}\begin{array}{rcl}
\Sha_{k,n-k}^{\phantom{\uparrow}}\Sigma_{k}^{\uparrow n-k}\Sha_{1,n-k-1}^{\phantom{\uparrow}}&=&
\Sha_{k,n-k}^{\phantom{\uparrow}}\Sha_{1,n-k-1}^{\phantom{\uparrow}}\Sigma_{k}^{\uparrow n-k}
\\[1em]
&=&\Sha_{k+1,n-k-1}^{\phantom{\uparrow}}\Sha_{k,1}^{\uparrow n-k-1}\Sigma_{k}^{\uparrow n-k}
=\Sha_{k+1,n-k-1}^{\phantom{\uparrow}}\Sigma_{k+1}^{\uparrow n-k-1}
\ .\end{array}\ee

\vskip .4cm
Shuffle elements find numerous applications in the theories of free Lie algebras, 
polylogarithms and multiple zeta values, Hopf algebras, differential calculus on quantum groups, 
homology of quantum Lie algebras, braidings of tensor spaces {\em etc.} 
\cite{R}, \cite{BB}, \cite{N}, \cite{Andr}, \cite{Rosso2}, \cite{Wor}, \cite{IO}, \cite{IOG}, 
\cite{GO}.

\section{B-shuffle algebras}

In this section we recall the definition of the Nichols--Woronowicz algebras and construct,
with the help of the baxterized elements, another family of graded associative algebras in the 
tensor spaces of local representations of the braid groups.

\paragraph{1.} Let $V$ be a vector space over a field ${\mathfrak k}$. For an operator 
$X\in {\mathrm{End}}(V^{\otimes j})$ we denote by the same symbol the operator 
$X\otimes {\mathrm{Id}}^{\otimes l}\in {\mathrm{End}}(V^{\otimes (j+l)})$ for any 
$l\in {\mathbb Z}_{\geq 0}$; $X^{\uparrow l}$ denotes the operator 
${\mathrm{Id}}^{\otimes l}\otimes X\in {\mathrm{End}}(V^{\otimes (j+l)})$.

\vskip .2cm
Let $\{ {\cal{T}}_{m,n}\}_{_{m,n\in {\mathbb Z}_{\geq 0}}}$ be a collection of operators, 
${\cal{T}}_{m,n}\in {\mathrm{End}}(V^{\otimes (m+n)})$, such that
\be\lb{ctshr} {\cal{T}}_{n+k,m}^{\phantom{\uparrow}}\, {\cal{T}}_{k,n}^{\uparrow m}=
{\cal{T}}_{k,m+n}^{\phantom{\uparrow}}\, {\cal{T}}_{n,m}^{\phantom{\uparrow}}
\qquad\forall\ m,n,k\in {\mathbb Z}_{\geq 0}\ .\ee
{}For tensors $u\in V^{\otimes m}$ and $v\in V^{\otimes n}$ let
\be\lb{shal} u\odot v:={\cal{T}}_{n,m}(u\otimes v)\in V^{\otimes (m+n)}\ .\ee
Due to (\ref{ctshr}), the space $\bigoplus_j V^{\otimes j}$ with the composition law $\odot$
is an associative graded algebra. Assume, in addition, that 
\be\lb{ico} {\cal{T}}_{m,0}={\mathrm{Id}}\quad {\mathrm{and}}\quad
{\cal{T}}_{0,m}={\mathrm{Id}} \qquad\forall\ m\in {\mathbb Z}_{\geq 0}\ee 
(then $1\in {\mathfrak k}\equiv V^{\otimes 0}$ is the identity element of the algebra). 
By (\ref{ctshr}) and (\ref{ico}), the following collection 
$\{ {\cal{S}}_m\}_{_{m\in {\mathbb Z}_{\geq 0}}}$ of operators,
${\cal{S}}_m\in {\mathrm{End}}(V^{\otimes m})$,
\be\lb{symshr} {\cal{S}}_0={\mathrm{Id}}\quad ,\quad {\cal{S}}_1={\mathrm{Id}}\quad ,\quad 
{\cal{S}}_{m+n}^{\phantom{\uparrow}}={\cal{T}}_{n,m}^{\phantom{\uparrow}}\,
{\cal{S}}_m^{\phantom{\uparrow}}\, {\cal{S}}_n^{\uparrow m}
\qquad\forall\ m,n\in {\mathbb Z}_{\geq 0}\ ,\ee
is well defined. The operation $\odot$ restricts on $\bigoplus_j {\mathrm{Im}}({\cal{S}}_j)$,
the direct sum of images of the operators ${\cal{S}}_j$, making it an associative graded algebra
as well.

\vskip .2cm
Let $\hat{R}\in {\mathrm{End}}(V\otimes V)$ be a solution of the Yang-Baxter equation, that is, 
$\hat{R}\hat{R}^{\uparrow 1}\hat{R}=\hat{R}^{\uparrow 1}\hat{R}\hat{R}^{\uparrow 1}$.
Denote by $\rho_{\hat{R}}$ the corresponding local representation of the braid groups $B_n$,
$\rho_{\hat{R}}(\sigma_i):=\hat{R}^{\uparrow (i-1)}$. Then the collection 
${\cal{T}}_{m,n}:=\rho_{\hat{R}}(\Sha_{m,n})$ obeys (\ref{ctshr}) and (\ref{ico}). The space
$\bigoplus_j {\mathrm{Im}}\rho_{\hat{R}}(\Sigma_j)$ with the composition law $\odot$ is called
Nichols--Woronowicz algebra.

\paragraph{2.}
The braid group rings admit a family of automorphisms
$\sigma_i\mapsto t\,\sigma_i$, where $t\in {\mathfrak k}^*$ is an arbitrary
parameter. The formal limits $\displaystyle{\lim_{t\rightarrow 0}}$ (the lowest power in $t$)
and $\displaystyle{\lim_{t\rightarrow \infty}}$ (the highest power in $t$) of the elements
$\Sigma_m$, $\Sha_{m,n}$ and the operation $\odot$ are well defined. For $t\rightarrow 0$
we obtain the usual tensor algebra, while for $t\rightarrow \infty$ the element 
$\Sigma_{n+1}$ becomes the lift of the longest element of the symmetric group 
${\mathbb{S}}_{n+1}$ to $B_{n+1}$,
\be\lb{ees4}\bar{\Sigma}_{n+1}=(\sigma_1\sigma_2\dots\sigma_n)(\sigma_1\dots\sigma_{n-1})
\dots (\sigma_1)\ .\ee
The shuffle elements, in the limit $\displaystyle{\lim_{t\rightarrow \infty}}$, become the
elements $\bar{\Sha}_{m,n}$ which, in a representation in a vector space $V$, equip
the tensor powers of $V$ with the standard braidings; the recurrency relations (\ref{shufrec}) 
and (\ref{shufrec2}) take the multiplicative form for $\bar{\Sha}_{m,n}$,
\be\lb{muh1}\bar{\Sha}_{m,n}=\bar{\Sha}_{m,n-1}\sigma_{m+n-1}\dots\sigma_n\ ,\ 
\bar{\Sha}_{m,n}=\bar{\Sha}_{m-1,n}^{\uparrow 1}\sigma_1\dots\sigma_n\ .\ee
Explicitly,
\be\bar{\Sha}_{m,n}=\left\{\begin{array}{l}\Bigl(\sigma_m\sigma_{m+1}\dots\sigma_{m+n-1}\Bigr)
\Bigl(\sigma_{m-1}\sigma_{m}\dots\sigma_{m+n-2}\Bigr)\Bigl(\sigma_1\sigma_2\dots\sigma_n
\Bigr)\ ,\\[1em]
\Bigl(\sigma_m\sigma_{m-1}\dots\sigma_1\Bigr)\Bigl(\sigma_{m+1}\sigma_m\dots\sigma_2\Bigr)
\Bigl(\sigma_{m+n-1}\sigma_{m+n-2}\dots\sigma_n\Bigr)\ .\end{array}\right.\ee
In addition to (\ref{muh1}), the elements $\bar{\Sha}_{m,n}$ satisfy
\be\lb{muh2}\bar{\Sha}_{m,n}=\sigma_m\dots\sigma_1\bar{\Sha}_{m,n-1}^{\uparrow 1}\ ,\ 
\bar{\Sha}_{m,n}=\sigma_m\dots\sigma_{m+n-1}\bar{\Sha}_{m-1,n}\ .\ee

\paragraph{3.}
In this section we shall construct another collection ${\cal{T}}_{m,n}$ starting with 
the elements $\sigma_k(x,y)$, satisfying the Yang-Baxter equation with spectral parameters 
\be\lb{sh1}\sigma_k(x_{k+1},x_{k+2})\sigma_{k+1}(x_k,x_{k+2})\sigma_k(x_k,x_{k+1})
=\sigma_{k+1}(x_k,x_{k+1})\sigma_k(x_k, x_{k+2})\sigma_{k+1}(x_{k+1},x_{k+2})\ee
and the locality condition
\be\lb{sh2}\sigma_k(x_{k},x_{k+1})\sigma_{l}(x_l,x_{l+1})=
\sigma_{l}(x_l,x_{l+1})\sigma_k(x_{k},x_{k+1})\qquad {\mathrm{if}}\ \ |k-l| > 1\ .\ee
Here $x_k$ are variables (spectral parameters). Depending on the situation, the elements 
$\sigma_k(x,y)$ can live in certain quotients of the braid group rings or be realized as operators.
We shall call $\sigma_k(x,y)$ baxterized elements (usually the term "baxterized" is applied
when $\sigma(x,y)$ is a function of the solution $\sigma$ of the constant Yang-Baxter 
equation).

\vskip .2cm
Let $\pi_k$ be the operator which permutes the variables $x_k$ and $x_{k+1}$,
$$\pi_kf(\dots,x_k,x_{k+1},\dots)=f(\dots,x_{k+1},x_k,\dots)\pi_k\ .$$
Relations (\ref{sh1}) and (\ref{sh2}) acquire the braid form (\ref{ryb1}) and (\ref{braidg})
for the elements
\be\lb{sh3}\underline{\sigma}_k:=\pi_k\sigma_k(x_k,x_{k+1})\ .\ee
The unitarity condition $\sigma_k(x_{k},x_{k+1})\sigma_k(x_{k+1},x_{k})=1$ (if imposed)
for the baxterized elements takes the form $\underline{\sigma}_k^2 =1$ for the elements 
(\ref{sh3}).

\vskip .2cm
The operators $\pi_k$ obey the braid group relations; prepare the elements 
$\bar{\Sha}_{m,n}\{\pi\}$ and $\bar{\Sigma}_m\{\pi\}$ from $\pi$'s; the elements 
$\bar{\Sha}_{m,n}\{\underline{\sigma}\}$ and $\bar{\Sigma}_m\{\underline{\sigma}\}$ built from 
$\underline{\sigma}$'s can be written, after moving all $\pi$'s to the left, in the form
\be\lb{preti}\bar{\Sha}_{m,n}\{\underline{\sigma}\} =
\bar{\Sha}_{m,n}\{\pi\}\widetilde{\Sha}_{m,n}(x_1,\dots,x_{m+n})\quad ,\quad
\bar{\Sigma}_m\{\underline{\sigma}\} =
\bar{\Sigma}_m\{\pi\}\widetilde{\Sigma}_m(x_1,\dots,x_m)\ ,\ee
where 
\be\lb{basha}\begin{array}{ccl}\widetilde{\Sha}_{m,n}(x_1,\dots,x_{m+n})&=&
\Bigl(\sigma_m(x_1,x_{m+n})\sigma_{m+1}(x_2,x_{m+n})
\dots\sigma_{m+n-1}(x_n,x_{m+n})\Bigr)\\[1em]
&\cdot&\Bigl(\sigma_{m-1}(x_1,x_{m+n-1})\sigma_{m}(x_2,x_{m+n-1})
\dots\sigma_{m+n-2}(x_n,x_{m+n-1})\Bigr)\\[1em]
&\dots&\Bigl(\sigma_1(x_1,x_{n+1})\sigma_2(x_2,x_{n+1})
\dots\sigma_n(x_n,x_{n+1})\Bigr)\end{array}\ee
and 
\be\lb{basi}\begin{array}{ccl}\widetilde{\Sigma}_m(x_1,\dots,x_m)&=&
\Bigl(\sigma_1(x_{m-1},x_m)\sigma_2(x_{m-2},x_m)\dots\sigma_{m-1}(x_1,x_m)\Bigr)\\[1em]
&\cdot&\Bigl(\sigma_1(x_{m-2},x_{m-1})\sigma_2(x_{m-3},x_{m-1})
\dots\sigma_{m-2}(x_1,x_{m-1})\Bigr)\dots\Bigl(\sigma_1(x_1,x_2)\Bigr)\ .\end{array}\ee

The elements $\bar{\Sha}_{m,n}\{\pi\}$ and $\bar{\Sigma}_n\{\pi\}$ are invertible and obey the 
relations (\ref{symsh}), (\ref{ctsh}), (\ref{muh1}) and (\ref{muh2}); substituting (\ref{preti}) 
into (\ref{symsh}), (\ref{ctsh}), (\ref{muh1}) and (\ref{muh2}), moving all $\pi$'s to the left
and simplifying, we find relations for $\widetilde{\Sigma}$'s and $\widetilde{\Sha}$'s alone.
The relations (\ref{muh1}) and (\ref{muh2}) take the form
\be\!\widetilde{\Sha}_{m,n}(x_1,\dots,x_{m+n})\!=\!\left\{\begin{array}{l}\!\!\!
\widetilde{\Sha}_{m,n-1}(\hat{x}_n)
\sigma_{m+n-1}(x_n,x_{m+n})\sigma_{m+n-2}(x_n,x_{m+n-1})\!\dots\!\sigma_n(x_n,x_{n+1}),\\[1em]
\!\!\!\widetilde{\Sha}_{m-1,n}^{\uparrow 1}(\hat{x}_{n+1})\sigma_1(x_1,x_{n+1})
\sigma_2(x_2,x_{n+1})\!\dots\!\sigma_n(x_n,x_{n+1}),\\[1em]
\!\!\!\sigma_m(x_1,x_{m+n})\sigma_{m-1}(x_1,x_{m+n-1})\!
\dots\!\sigma_1(x_1,x_{n+1})\widetilde{\Sha}_{m,n-1}^{\uparrow 1}(\hat{x}_1),\\[1em]
\!\!\!\sigma_m(x_1,x_{m+n})\sigma_{m+1}(x_2,x_{m+n})\!\dots\!\sigma_{m+n-1}(x_n,x_{m+n})
\widetilde{\Sha}_{m-1,n}(\hat{x}_{m+n}\}),\end{array}\right.\ee
where "$\hat{x}_j$" means that the argument $x_j$ is omitted. For a set 
$\overrightarrow{x}=\{ x_1,\dots,x_n\}$ of arguments, let $\overleftarrow{x}:=\{ x_n,\dots,x_1\}$
be the reversed set. The relation (\ref{symsh}) becomes:
\be\lb{wtsi}\widetilde{\Sigma}_{m+n}^{\phantom{\uparrow}}(\overrightarrow{x},\overrightarrow{y})
=\widetilde{\Sha}_{n,m}^{\phantom{\uparrow}}(\overleftarrow{x},\overleftarrow{y})\,
\widetilde{\Sigma}_m^{\phantom{\uparrow}}(\overrightarrow{x})\,
\widetilde{\Sigma}_n^{\uparrow a}(\overrightarrow{y})\ ,\ee
where $\overrightarrow{x}=\{ x_1,\dots,x_m\}$ and $\overrightarrow{y}=\{ y_1,\dots,y_n\}$;
the relation (\ref{ctsh}) becomes
\be\widetilde{\Sha}_{n+k,m}^{\phantom{\uparrow}}
(\overrightarrow{x},\overrightarrow{z},\overrightarrow{y})
\widetilde{\Sha}_{k,n}^{\uparrow m}(\overrightarrow{y},\overrightarrow{z})=
\widetilde{\Sha}_{k,m+n}^{\phantom{\uparrow}}
(\overrightarrow{y},\overrightarrow{x},\overrightarrow{z})
\widetilde{\Sha}_{n,m}^{\phantom{\uparrow}}(\overrightarrow{x},\overrightarrow{y})\ ,\ee
where $\overrightarrow{x}=\{ x_1,\dots,x_m\}$, $\overrightarrow{y}=\{ y_1,\dots,y_n\}$
and $\overrightarrow{z}=\{ z_1,\dots,z_k\}$.

\vskip .2cm
After the removal of all $\pi$'s, one can give values to the spectral variables. Each
$\widetilde{\Sigma}_m$ can be evaluated on its own sequence 
$\vec{x}^{(m)}=(x_1^{(m)},\dots,x_m^{(m)})$. In the relation (\ref{wtsi}), the 
beginning of the sequence for $\widetilde{\Sigma}_{m+n}$ becomes the beginning of
the sequence for $\widetilde{\Sigma}_m$ while its end becomes the beginning of the
sequence for $\widetilde{\Sigma}_n$. This is a strong restriction; if it is imposed on the
sequences themselves, the general solution is that each $x_j^{(m)}$ is equal to one and the same 
number. However, assume that the baxterization is "trigonometric", the baxterized elements depend 
on the ratio of the spectral parameters, $\sigma(x,y)=\sigma(x/y)$.  The Yang-Baxter equation then 
reads
\be\lb{ybeH}\sigma_n(x)\sigma_{n-1}(xy)\sigma_n(y)=\sigma_{n-1}(y)\sigma_n(xy)\sigma_{n-1}(x)\ .
\ee
Now $\widetilde{\Sigma}_m(\vec{x}^{(m)})=\widetilde{\Sigma}_m(\alpha\vec{x}^{(m)})$ for an arbitrary constant $\alpha\neq 0$ and the general solution of the restrictions imposed
by (\ref{wtsi}) for the projectivized sequences is $(x_1^{(m)},\dots,x_m^{(m)})=(1,s^{-1},s^{-2},\dots,s^{1-m})$, the 
geometric progression. Denote $\widetilde{\Sigma}_m(1,s^{-1},s^{-2},\dots,s^{1-m})$ by
$^{^s}\!\Sigma_m$ and 
$\widetilde{\Sha}_{m,n}(s^{1-n},\dots,1,s^{1-m-n},\dots,s^{-n})$ by 
$^{^s}\!\Sha_{m,n}$. Explicitly:
\be\lb{basi2} ^{^s}\!\Sigma_m=
\Bigl(\sigma_1(s)\sigma_2(s^2)\dots\sigma_{m-1}(s^{m-1})\Bigr)
\Bigl(\sigma_1(s)\sigma_2(s^2)\dots\sigma_{m-2}(s^{m-2})\Bigr)\dots\Bigl(\sigma_1(s)\Bigr)
\ee
and
\be\lb{basha2} ^{^s}\!\Sha_{m,n}=\Bigl(\sigma_m(s)\dots\sigma_{m+n-1}(s^n)\Bigr)
\Bigl(\sigma_{m-1}(s^2)\dots\sigma_{m+n-2}(s^{n+1})\Bigr)
\Bigl(\sigma_1(s^m)\dots\sigma_n(s^{m+n-1})\Bigr)\ .\ee

The elements $^{^s}\!\Sha_{m,n}$ obey the relation (\ref{ctsh}). Therefore, in a local
representation $\rho$, the collection ${\cal{T}}_{m,n}:=\rho_{\hat{R}}(^{^s}\!\Sha_{m,n})$ 
obeys (\ref{ctshr}) and (\ref{ico}) and defines a one parameter family of graded associative 
algebras on $\bigoplus_j V^{\otimes j}$ together with the subalgebras on 
$\bigoplus_j {\mathrm{Im}}({\cal{S}}_j)$ (now ${\cal{S}}_m=\rho_{\hat{R}}(^{^s}\!\Sigma_m)$), 
which we propose to call {\em b-shuffle algebras} ("b" from "baxterized"; maybe the 
term "buffle" would be an apt acronym). 

\vskip .2cm
It is known that the element $\bar{\Sigma}$ admits reduced expressions starting (or ending) 
with $\sigma_j$ for every $j=1,\dots,m-1$. In particular, $\bar{\Sigma}\{\underline{\sigma}\}$ 
can start (or end) with every $\underline{\sigma}_j$. It follows that  
$\widetilde{\Sigma}_m(x_1,\dots,x_m)$ can start (or end) with $\sigma_j(x_j,x_{j+1})$ for
every $j$. We shall use this for the trigonometric $\sigma$'s:
\be\lb{retr} ^{^s}\!\Sigma_m\ \ {\mathrm{has\ a\ reduced\ expression\ of\ the\ form}}\ \sigma_j(s)\cdot(\dots)
\ {\mathrm{or}}\ (\dots)\cdot\sigma_j(s)\ \ \forall\ j=1,\dots,m-1\ .\ee

\vskip .2cm
The baxterization is known for the Hecke and BMW quotients of the braid group rings;
it is trigonometric. In the next section we discuss the baxterized collections for these quotients.

\paragraph {4. Remarks.} {\bf (a)} We suggest another natural source for collections ${\cal{T}}_{m,n}$ satisfying (\ref{ctshr}) and (\ref{ico}). 

\vskip .2cm
Let ${\cal{A}}$ be a Hopf algebra. Assume that ${\cal{A}}$ admits a twist ${\cal{F}}$, that is, 
an element ${\cal{F}}\in {\cal{A}}\otimes {\cal{A}}$ which satisfies 
\be\lb{bai} {\cal{F}}\cdot (\Delta\otimes {\mathrm{Id}})({\cal{F}})=
{\cal{F}}^{\uparrow 1}\cdot ({\mathrm{Id}}\otimes\Delta)({\cal{F}})\ .\ee 
Here $\Delta$ is the coproduct and $^\uparrow$ is the shift in the copies of ${\cal{A}}$ in 
${\cal{A}}^{\otimes j}$. Define ${\cal{F}}_{m,0}:=1$, ${\cal{F}}_{0,m}:=1$,
$m\in {\mathbb Z}_{\geq 0}$, and
\be {\cal{F}}_{m,n}:=\Delta^{m-1}\otimes\Delta^{n-1}({\cal{F}})\qquad ,\ \  
m,n\in {\mathbb Z}_{\geq 1}\ .\ee
It is straightforward to verify that 
\be\lb{utw}{\cal{F}}_{k,m}^{\phantom{\uparrow}}{\cal{F}}_{m+k,n}^{\phantom{\uparrow}}
={\cal{F}}_{m,n}^{\uparrow k}{\cal{F}}_{k,m+n}^{\phantom{\uparrow}}\ee
(for $m=n=k=1$ this is (\ref{bai}); by induction, ${\mathrm{Id}}^{i-1}\otimes\Delta
\otimes {\mathrm{Id}}^{m+n+k-i}$ increases $k$ by 1 for $1\leq i\leq k$, 
$m$ by 1 for $k<i\leq m+k$ and $n$ by 1 for $m+k<i\leq m+n+k$).

\vskip .2cm
Therefore, given a representation $\rho$ of ${\cal{A}}$, the relations (\ref{ctshr}) and 
(\ref{ico}) hold for 
$${\cal{T}}_{m,n}:=\varpi\circ\rho^{\otimes (m+n)}({\cal{F}}_{n,m})\ ,$$
where $\varpi$ is any operation which reverses the order of terms in both sides of (\ref{utw})
(it can be a transposition or, if ${\cal{F}}_{m,n}$ are invertible for all $m$ and $n$, an inversion).

\vskip .2cm
It might be of interest to investigate this type of collections ${\cal{T}}_{m,n}$
for the twists \cite{IOt} corresponding to Belavin--Drinfeld triples.

\paragraph{(b)} Assume, in addition, that ${\cal{F}}$ satisfies
\be (\Delta\otimes {\mathrm{Id}})({\cal{F}})={\cal{F}}_{\{1,3\}}{\cal{F}}_{\{2,3\}}\ ,\ 
({\mathrm{Id}}\otimes\Delta)({\cal{F}})={\cal{F}}_{\{1,3\}}{\cal{F}}_{\{1,2\}}\ ,\ee
where ${\cal{F}}_{\{i,j\}}$ is the element ${\cal{F}}$ located in the copies number $i$ and $j$
in ${\cal{A}}\otimes {\cal{A}}\otimes\dots$ ; for example, for a quasi-triangular Hopf algebra, 
${\cal{F}}$ can be the universal $R$-matrix. Then
\be\lb{fumn} {\cal{F}}_{m,n}=\Bigl({\cal{F}}_{\{1,m+n\}}\dots {\cal{F}}_{\{m,m+n\}}
\Bigr)\Bigl({\cal{F}}_{\{1,m+n-1\}}\dots {\cal{F}}_{\{m,m+n-1\}}\Bigr)\dots
\Bigl({\cal{F}}_{\{1,m+1\}}\dots {\cal{F}}_{\{m,m+1\}}\Bigr)\ee
(in each bracket the first index increases from 1 to $m$, the second one is constant); this 
formula generalizes the formula $\Delta\otimes\Delta ({\cal{R}})={\cal{R}}_{\{1,4\}}
{\cal{R}}_{\{2,4\}}{\cal{R}}_{\{1,3\}}{\cal{R}}_{\{2,3\}}$ used in the theory of 
quasi-triangular Hopf algebras for establishing properties of the element giving the square of the
antipode by conjugation, see, e.g., \cite{O}, chapter 4. It follows from (\ref{fumn}) that
\be\lb{inuf}{\cal{F}}_{m,n}=\Bigl({\cal{F}}_{\{1,m+n\}} {\cal{F}}_{\{1,m+n-1\}}\dots 
{\cal{F}}_{\{1,m+1\}}\Bigr) {\cal{F}}_{m-1,n}^{\uparrow 1}\ .\ee
Given a representation $\rho$ of ${\cal{A}}$ on a vector space $V$, let $P_i$ be the flip operator
in the copies number $i$ and $i+1$ of the space $V$ in $V\otimes V\otimes\dots$ ; let 
$F:=\rho^{\otimes 2}({\cal{F}})$ and $\hat{F}:=P_1F$; for an operator 
$X\in {\mathrm{End}}(V\otimes V)$ denote by $X_{\{i,j\}}$ the operator $X$ acting in the
copies number $i$ and $j$ of the space $V$ in $V\otimes V\otimes\dots$ and let 
$X_i:=X_{\{i,i+1\}}$. Then
\be\lb{unf}\rho^{\otimes (m+n)}({\cal{F}}_{m,n})=
\bar{\Sha}_{m,n}\{ P\}\bar{\Sha}_{n,m}\{\hat{F}\}\ ,\ee
where $\bar{\Sha}_{m,n}\{ P\}$ are built from $P$'s and $\bar{\Sha}_{n,m}\{\hat{F}\}$ from 
$\hat{F}$'s. Indeed, by (\ref{inuf}) and induction,
\be\begin{array}{l} \rho^{\otimes (m+n)}({\cal{F}}_{m,n})=\Bigl(F_{\{1,m+n\}}\dots F_{\{1,m+1\}}
\Bigr)\bar{\Sha}_{m-1,n}^{\uparrow 1}\{ P\}\bar{\Sha}_{n,m-1}^{\uparrow 1}\{\hat{F}\}\\[1em]
\hspace{1cm}\! =\!\Bigl(P_{\{1,m+n\}}\dots P_{\{1,m+1\}}\Bigr)\!\Bigl({\hat{F}}_{m+n-1} 
{\hat{F}}_{m+n-2}\dots {\hat{F}}_{m+1}{\hat{F}}_{\{1,m+1\}}\Bigr)\!
\bar{\Sha}_{m-1,n}^{\uparrow 1}\{ P\}\bar{\Sha}_{n,m-1}^{\uparrow 1}\{\hat{F}\}\ .\end{array}\ee
Use now 
\be \Bigl({\hat{F}}_{m+n-1} {\hat{F}}_{m+n-2}\dots  
{\hat{F}}_{m+1}{\hat{F}}_{\{1,m+1\}}\Bigr)\bar{\Sha}_{m-1,n}^{\uparrow 1}\{ P\}=
\bar{\Sha}_{m-1,n}^{\uparrow 1}\{ P\}\Bigl({\hat{F}}_n {\hat{F}}_{n-1}\dots  
{\hat{F}}_{1}\Bigr)\ ,\ee
the first recursion relelation in (\ref{muh2}) for $\bar{\Sha}_{n,m}\{\hat{F}\}$ and
\be\begin{array}{l} \Bigl(P_{\{1,m+n\}}\dots P_{\{1,m+1\}}\Bigr)
\bar{\Sha}_{m-1,n}^{\uparrow 1}\{ P\} =\bar{\Sha}_{m-1,n}^{\uparrow 1}\{ P\}
\Bigl(P_{\{1,n+1\}}P_{\{1,n\}}\dots P_{\{1,2\}}\Bigr)\\[1em]
\hspace{2cm}=\bar{\Sha}_{m-1,n}^{\uparrow 1}\{ P\}\Bigl(P_1P_2\dots P_n\Bigr)
=\bar{\Sha}_{m,n}\{ P\}\end{array}\ee
(by the second recurrency relation in (\ref{muh1}) for $\bar{\Sha}_{m,n}\{ P\}$) to finish the 
proof of (\ref{unf}).

\vskip .2cm
Thus the elements ${\cal{F}}_{m,n}$ can be regarded as the universal (in the Hopf
algebra theoretical sense) counterpart of the elements $\bar{\Sha}_{m,n}$.

\paragraph{(c)} We describe an operation which transforms a collection ${\cal{T}}_{m,n}$ 
satisfying (\ref{ctshr}) and (\ref{ico}) into another, "dual", collection $\check{\cal{T}}_{m,n}$
satisfying (\ref{ctshr}) and (\ref{ico}). 

\vskip .2cm
Keep the notation from the previous remark. Let $X\in {\mathrm{End}}(V^{\otimes m})$ and $Y\in {\mathrm{End}}(V^{\otimes n})$ be two 
operators. Then
\be\lb{pesha}\bar{\Sha}_{m,n}\{ P\}X^{\uparrow n}Y=XY^{\uparrow m}\bar{\Sha}_{m,n}\{ P\}\ .\ee
Define $\check{\cal{T}}_{m,n}$ by
\be\lb{dedu} {\cal{T}}_{m,n}:=\check{\cal{T}}_{n,m}\bar{\Sha}_{m,n}\{ P\}\ \ {\mathrm{or}}\ \ 
\check{\cal{T}}_{m,n}:={\cal{T}}_{n,m}\bar{\Sha}_{m,n}\{ P\}\ .\ee
The equivalence of two definitions follows from
\be \bar{\Sha}_{m,n}\{ P\}^{-1}=\bar{\Sha}_{n,m}\{ P\}\ .\ee

The relation (\ref{ico}) is satisfied for the collection $\check{\cal{T}}_{m,n}$. The
relation (\ref{ctshr}) reads, by (\ref{pesha}),
\be \check{\cal{T}}_{m,n+k}\check{\cal{T}}_{n,k}\bar{\Sha}_{n+k,m}\{ P\}
\bar{\Sha}_{k,n}^{\uparrow m}\{ P\}=\check{\cal{T}}_{m+n,k}
\check{\cal{T}}_{m,n}^{\uparrow k}\bar{\Sha}_{k,m+n}\{ P\}\bar{\Sha}_{n,m}\{ P\}\ .\ee
Since
\be\bar{\Sha}_{n+k,m}\{ P\}\bar{\Sha}_{k,n}^{\uparrow m}\{ P\}
=\bar{\Sha}_{k,m+n}\{ P\}\bar{\Sha}_{n,m}\{ P\}\ee
it follows that the relation (\ref{ctshr}) is as well satisfied for the collection 
$\check{\cal{T}}_{m,n}$.

\vskip .2cm
With the help of the identity $ \bar{\Sigma}_m\{ P\}^2={\mathrm{Id}}$, it is straightforward 
to verify that the collection $\check{\cal{S}}_m$ for $\check{\cal{T}}_{m,n}$ is given by
\be \check{\cal{S}}_m={\cal{S}}_m\bar{\Sigma}_m\{ P\}\ .\ee

\section{Hecke and BMW algebras\vspace{.25cm}}

\vspace{-.2cm}
In the sequel we call the elements $\Sha_{m,n}$ additive shuffles and $^{^s}\!\Sha_{m,n}$
multiplicative shuffles. In this section we derive the sequences of the (anti-)symmetrizers 
for the Hecke and BMW algebras with the help of the multiplicative shuffles. We compare the
multiplicative versions with known expressions for the (anti-)symmetrizers. 

\vskip .2cm
We derive a new expression for the (anti-)symmetrizers in terms of the highest multiplicative 
1-shuffles alone and, for the Hecke algebras, in terms of the highest additive 1-shuffles alone.

\vskip .2cm
In principle, the Hecke algebras can be considered as quotients of the BMW algebras
and many formulas for the Hecke algebras can be obtained from this point of view. 
Because of importance of the Hecke algebras we prefer however to treat them separately.

\subsection{Hecke algebras}

\paragraph{1.} The tower of the $A$-Type Hecke algebras $H_{M+1}(q)$ (see {\em e.g.} \cite{Jon1} 
and references therein) depends on a parameter $q\in {\mathfrak{k}}^*$; the algebra  $H_{M+1}(q)$ 
is the quotient of the braid group ring ${\mathfrak{k}}B_{M+1}$ by
\be\lb{ahecke}\sigma^2_i=(q-q^{-1})\sigma_i+1\ \ ,\ i=1,\dots,M\ .\ee

{}For $q^2\neq 1$, the baxterized elements have the form
\be\lb{baxtH}\sigma_i(x):=\frac{1}{q-q^{-1}}\,(x\sigma_i-x^{-1}\sigma_i^{-1})\ ;\ee
they are normalized, $\sigma_i(1)=1$, and satisfy  the unitarity condition
\be\lb{hunit}\sigma_i(x)\sigma_i(x^{-1})=1-\frac{(x-x^{-1})^2}{(q-q^{-1})^2}\ .\ee

\paragraph{2.} The symmetrizers $S_n$, $n=1,\dots,M+1$, (\cite{Jon0}, \cite{W1}, \cite{Gur}) are 
the non-zero elements, 
\be\lb{desy} S_1=1\ ,\ S_n\in H_n(q)\subset H_{M+1}(q)\ ,\ee
which satisfy 
\be\lb{desy'} \sigma_iS_n=S_n\sigma_i=qS_n\ ,\ i=1,\dots,n-1\ ,\ee
which forces $S_n^2\sim S_n$; they are normalized by
\be\lb{desy2} S_n^2=S_n\ .\ee
The sequence $\{ S_n\}$ is defined by (\ref{desy}), (\ref{desy'}) and (\ref{desy2}) uniquely 
(and it does exist for the Hecke quotients for generic $q$); the anti-symmetrizers are 
related to the symmetrizers by the isomorphisms $H_{M+1}(q)\rightarrow H_{M+1}(-q^{-1})$, 
$\sigma_i\mapsto\sigma_i$. 

\paragraph{3.} The symmetrizers can be quickly constructed with the help of the baxterized
elements. Let $[n]_q:=(q^n-q^{-n})/(q-q^{-1})$, $[n]_q!:=[1]_q[2]_q\cdots [n]_q$ and
$[n]_q^\$ :=[1]_q![2]_q!\cdots [n]_q!$. By (\ref{retr}), 
$\sigma_i\, ^{^q}\!\Sigma_n=\, ^{^q}\!\Sigma_n\sigma_i=q\, ^{^q}\!\Sigma_n$, or, 
equivalently, $\sigma_i(q) ^{^q}\!\Sigma_n=[i+1]_q\, ^{^q}\!\Sigma_n$, so
$(^{^q}\!\Sigma_n)^2=[n]_q^\$\, ^{^q}\!\Sigma_n$ and
\be\lb{mshu} S_n=\frac{1}{[n]_q^\$}\  ^{^q}\!\Sigma_n\ .\ee
In particular, the symmetrizers satisfy the recurrent relation
\be\lb{recm} S_n=\frac{1}{[n]_q!}\  ^{^q}\!\Sha_{1,n-1}S_{n-1}\ .\ee

\paragraph{4.}
We recall several other forms of the symmetrizers and compare them with (\ref{mshu}) and
(\ref{recm}). A convenient recurrent relation for the symmetrizers is (see e.g. \cite{Jimb1}, 
\cite{HIOPT}):
\be\lb{santis2}S_n=S_{n-1}\,\frac{\sigma_{n-1}(q^{n-1})}{[n]_q}\, S_{n-1}\ .\ee
This is checked either by verifying (\ref{desy}), (\ref{desy'}) and (\ref{desy2}) and then 
by uniqueness or, using (\ref{recm}), by the following calculation:
\be\lb{dersan}\begin{array}{rcl} [n]_q!S_n&=&\sigma_1(q)\dots\sigma_{n-2}(q^{n-2})
\sigma_{n-1}(q^{n-1})S_{n-1}=\sigma_1(q)\dots\sigma_{n-2}(q^{n-2})\sigma_{n-1}(q^{n-1})
S_{n-2}S_{n-1}\\[1em]
&=&\sigma_1(q)\dots\sigma_{n-2}(q^{n-2})S_{n-2}\sigma_{n-1}(q^{n-1})S_{n-1}=
S_{n-1}\sigma_{n-1}(q^{n-1})S_{n-1}\ .\end{array}\ee

Denote $\Sha_{1,n}\{ q\sigma\}$ (the additive shuffle built with the $q\sigma_1,\dots,
q\sigma_{n-1}$) by $\sha_{1,n}$. There is another recurrent relation for the symmetrizers in 
terms of the additive shuffles:
\be\lb{reca} S_n=\frac{q^{1-n}}{[n]_q}\ \sha_{1,n-1}S_{n-1}\ .\ee
In other words,
\be\lb{santis} S_n=\frac{q^{-\frac{n(n-1)}{2}}}{[n]_q!}\ \Sigma_n\{ q\sigma\}\ .\ee
This is checked again either by verifying (\ref{desy}), (\ref{desy'}) and (\ref{desy2}) and 
then by uniqueness or, using (\ref{santis2}), by induction:
\be\lb{santider}\begin{array}{rcl} [n]_qS_n&=&S_{n-1}\sigma_{n-1}(q^{n-1})S_{n-1}=
\frac{q^{2-n}}{[n-1]_q}\sha_{1,n-2} S_{n-2}\sigma_{n-1}(q^{n-1})S_{n-1}\\[1em]
&=&\frac{q^{2-n}}{[n-1]_q}\sha_{1,n-2}\sigma_{n-1}(q^{n-1})S_{n-1}
=\frac{q^{2-n}}{[n-1]_q}\sha_{1,n-2} \bigl( [n-1]_q\sigma_{n-1}+q^{1-n}\bigr) S_{n-1}\\[1em]
&=&\frac{q^{2-n}}{[n-1]_q}\bigl( [n-1]_q\sha_{1,n-2}\sigma_{n-1}+
q^{-1}[n-1]_q\bigr) S_{n-1}=q^{1-n}\sha_{1,n-1}S_{n-1}\ .\end{array}\ee
We stress that the factors $\frac{1}{[n]_q!}\  ^{^q}\!\Sha_{1,n-1}$ in (\ref{recm}) and 
$\frac{q^{1-n}}{[n]_q}\ \sha_{1,n-1}$ in (\ref{reca}) differ; the multiplicative and
additive shuffles do not coincide although their products -- the symmetrizers -- do.

\paragraph{5.}It turns out that the symmetrizer $S_n$ can be expressed in terms of the 
multiplicative 1-shuffle $^{^q}\!\Sha_{1,n-1}$ or in terms of the additive 1-shuffle 
$\sha_{1,n-1}$ only.

\vskip .2cm
{}For the additive shuffle, we prove by induction that, for $k=1,\dots,n-1$,
\be\lb{adle}\prod_{j=0}^{k-1}(\sha_{1,n-1}-q^{j-1}[j]_q)
=q^{k(k-1)}\sha_{1,n-1}\sha_{1,n-2}\dots\sha_{1,n-k}\ .\ee
{}For $k=1$ there is nothing to prove. Assume that (\ref{adle}) holds for some $k<n-1$. 
By (\ref{lue}), the right hand side is divisible, from the right, by $S_k^{\uparrow (n-k)}$.
Multiply (\ref{adle}) by the factor $(\sha_{1,n-1}-q^{k-1}[k]_q)$ from the right and substitute,
in the right hand side,
$$\sha_{1,n-1}=\sha_{1,k-1}^{\uparrow (n-k)}+q^k\sha_{1,n-1-k}\sigma_{n-k}\dots\sigma_{n-1}$$
in this factor. Since $S_k^{\uparrow (n-k)}\sha_{1,k-1}^{\uparrow (n-k)}=q^{k-1}[k]_q
S_k^{\uparrow (n-k)}$, the product in the right hand side simplifies,
\be\begin{array}{c}\sha_{1,n-1}\sha_{1,n-2}\dots\sha_{1,n-k}(-q^{k-1}[k]_q+
\sha_{1,k-1}^{\uparrow (n-k)}+q^k\sha_{1,n-1-k}\sigma_{n-k}\dots\sigma_{n-1})\\[1em]
=q^k\sha_{1,n-1}\sha_{1,n-2}\dots\sha_{1,n-k}\sha_{1,n-1-k}\sigma_{n-k}\dots\sigma_{n-1}
=q^{2k}\sha_{1,n-1}\sha_{1,n-2}\dots\sha_{1,n-1-k}\end{array}\ee
(in the last equality we again used (\ref{lue}) for 
$\sha_{1,n-1}\sha_{1,n-2}\dots\sha_{1,n-k}\sha_{1,n-1-k}$), establishing the induction step.

\vskip .2cm
In particular, at $k=n-1$, we obtain, by (\ref{reca}), the expression of $S_n$ in terms
of $\sha_{1,n-1}$,
\be\lb{sysha} S_n=\frac{q^{-\frac{(n-1)(3n-4)}{2}}}{[n]_q!}\ 
\prod_{j=0}^{n-2}(\sha_{1,n-1}-q^{j-1}[j]_q)\ .\ee

\paragraph{6.} For the multiplicative shuffle, we prove by induction that, for $k=1,\dots,n$,
\be (\, ^{^q}\!\Sha_{1,n})^k=\frac{[n+1-k]_q^\$\ ([n+1]_q!)^k}{[n]_q^\$}\, ^{^q}\!\Sha_{1,n}
\, ^{^q}\!\Sha_{1,n-1}\dots\, ^{^q}\!\Sha_{1,n+1-k}\ .\ee
{}For $k=1$ there is nothing to prove. Assume that (\ref{adle}) holds for some $k<n$.
The relations (\ref{symsh}) and (\ref{ctsh}) hold for $\Sigma_m=\, ^{^q}\!\Sigma_m$
and $\Sha_{m,n}=\, ^{^q}\!\Sha_{m,n}$. Therefore, (\ref{lue}) holds as well and
the product $\, ^{^q}\!\Sha_{1,n}\, ^{^q}\!\Sha_{1,n-1}\dots\, ^{^q}\!\Sha_{1,n-k}$
is divisible, from the right, by $\, ^{^q}\!\Sigma_{k+1}^{\uparrow (n-k)}$. The induction step is:
\be\begin{array}{rcl}\, ^{^q}\!\Sha_{1,n}\dots\, ^{^q}\!\Sha_{1,n+1-k}\cdot
\, ^{^q}\!\Sha_{1,n}&=&\, ^{^q}\!\Sha_{1,n}\dots\, ^{^q}\!\Sha_{1,n+1-k}
\, ^{^q}\!\Sha_{1,n-k}\cdot\sigma_{n-k+1}(q^{n-k+1})\dots\sigma_n(q^n)\\[1em]
&=&{\displaystyle\frac{[n+1]_q!}{[n+1-k]_q!}}\, ^{^q}\!\Sha_{1,n}\, ^{^q}\!\Sha_{1,n-1}\dots
\, ^{^q}\!\Sha_{1,n+1-k}\, ^{^q}\!\Sha_{1,n-k}\ ,\end{array}\ee
since $S_{n+1}\sigma_k(q^k)=[k+1]_qS_{n+1}$, $k=1,\dots,n$.

\vskip .2cm
In particular, at $k=n$, we obtain, by (\ref{mshu}), the expression of $S_n$ in terms
of $\, ^{^q}\!\Sha_{1,n}$,
\be\lb{syshh}S_{n+1}=\Bigl(\frac{1}{[n+1]_q!}\, ^{^q}\!\Sha_{1,n}\Bigr)^n\ .\ee

\paragraph {7. Remark.} Let $\hat{R}$ be a Hecke Yang--Baxter matrix and $\rho_{\hat{R}}$ 
the corresponding local representation of the tower of the Hecke algebras. The symmetrizers
${\cal{S}}_j$ built with ${\cal{T}}_{m,n}'=\rho_{\hat{R}}(\, ^{^s}\!\Sha_{m,n})$, at
$s=q$, are the same as the symmetrizers built with 
${\cal{T}}_{m,n}''=\rho_{\hat{R}}(\Sha_{m,n}\{ t\sigma\} )$, at $t=q$ (the symmetrizers
coincide at $s^2=q^2$ and $t=q$ or $s^2=q^{-2}$ and $t=-q^{-1}$, these are the values for 
the anti-symmetrizers; the symmetrizers coincide trivially at $s^2=1$ and $t=0$; otherwise 
the symmetrizers for $\{{\cal{T}}_{m,n}'\}$ and $\{{\cal{T}}_{m,n}''\}$ differ). Therefore, 
for the Hecke algebras, the b-shuffle algebra on $\bigoplus_j {\mathrm{Im}}({\cal{S}}_j)$ 
coincides with the Nichols--Woronowicz algebra (or the symmetric algebra of the quantum space). 
Indeed, the composition law (\ref{shal}) on $\bigoplus_j {\mathrm{Im}}({\cal{S}}_j)$ can be 
written in the following equivalent form:
\be\lb{fe1}u\odot v:={\cal{S}}_{m+n}(u'\otimes v')\ ,\ee
where $u={\cal{S}}_mu'$ and $v={\cal{S}}_nv'$. Also,
${\mathrm{Im}}({\cal{S}}_j)\simeq V^{\otimes j}/{\mathrm{Ker}}({\cal{S}}_j)$, and
the algebra on $\bigoplus_j {\mathrm{Im}}({\cal{S}}_j)$ can be defined alternatively as
the algebra on $\bigoplus_j V^{\otimes j}/{\mathrm{Ker}}({\cal{S}}_j)$ with the composition
law
\be\lb{fe2}\bar{u}\circ \bar{v}:=u\otimes v\ \ {\mathrm{mod}}
\ {\mathrm{Ker}}({\cal{S}}_{m+n})\ ,\ee
where $\bar{u}\in V^{\otimes m}/{\mathrm{Ker}}({\cal{S}}_m)$ and
$\bar{v}\in V^{\otimes n}/{\mathrm{Ker}}({\cal{S}}_n)$; $u\in V^{\otimes m}$ and 
$v\in V^{\otimes n}$ are representatives of $\bar{u}$ and $\bar{v}$, respectively.
In the formulations (\ref{fe1}) or (\ref{fe2}), the algebra on 
$\bigoplus_j {\mathrm{Im}}({\cal{S}}_j)$ or
$\bigoplus_j V^{\otimes j}/{\mathrm{Ker}}({\cal{S}}_j)$
depends only on the collection $\{ {\cal{S}}_j\}$; the composition laws (\ref{fe1}) or 
(\ref{fe2}) are well defined if ${\cal{S}}_{m+n}$ is divisible by ${\cal{S}}_m$ and
${\cal{S}}_n^{\uparrow m}$ from the right (which is, in general, weaker than
${\cal{S}}_{m+n}^{\phantom{\uparrow}}={\cal{T}}_{n,m}^{\phantom{\uparrow}}\,
{\cal{S}}_m^{\phantom{\uparrow}}\, {\cal{S}}_n^{\uparrow m}$); when, say, the representative 
$u$ of $\bar{u}$ changes, $u\sim u+\delta u$, ${\cal{S}}_m(\delta u)=0$, so 
$\delta u\otimes v\in {\mathrm{Ker}}({\cal{S}}_{m+n})$ and the product $\bar{u}\circ \bar{v}$
does not change, $u\otimes v\equiv (u+\delta u)\otimes v\ \ {\mathrm{mod}}
\ {\mathrm{Ker}}({\cal{S}}_{m+n})$.

\vskip .2cm
However, the algebras on the space $\bigoplus_j V^{\otimes j}$, built with ${\cal{T}}_{m,n}'$
or ${\cal{T}}_{m,n}''$, are very different, as it is seen, for example from the comparison 
of the spectra of the multiplicative and additive 1-shuffles in Section \ref{secspec}.
The collections $\{ {\cal{T}}_{m,n}'\}$ and $\{ {\cal{T}}_{m,n}''\}$ seem to have
different ranges of applicability (already for the BMW algebras, the symmetrizers for these
two collections do not coincide). 

\subsection{BMW algebras}

The tower of the Birman-Murakami-Wenzl algebras $B\!M\!W_{M+1}(q,\nu)$ was introduced 
in \cite{M1} and \cite{BW}; it depends on two parameters, $q\in {\mathfrak{k}}^*$ and
$\nu\in {\mathfrak{k}}\setminus \{0,q,-q^{-1}\}$. For $q^2\neq 1$, the algebra 
$B\!M\!W_{M+1}(q,\nu)$ is the quotient of the braid group ring ${\mathfrak{k}}B_{M+1}$ by
\be\lb{bmw1}\kappa_i\sigma_i=\sigma_i\kappa_i=\nu\kappa_i\ ,\ee
\be\lb{bmw2}\kappa_i\sigma_{i-1}\kappa_i=\nu^{-1}\kappa_i\quad ,\quad 
\kappa_i\sigma_{i-1}^{-1}\kappa_i=\nu\kappa_i\ ,\ee
where the elements $\kappa_i$ are defined by
\be\lb{bmw3}\sigma_i-\sigma_i^{-1}=(q-q^{-1})(1-\kappa_i)\ .\ee

The Hecke quotient is $\kappa_i=0$.

\vskip .2cm
{}For $q^2\neq 1$, the baxterized elements have the form
(\cite{Jon2}, \cite{Mu}, \cite{J}, \cite{Is1})
\be\lb{bmwbax}\sigma_i(x):=x^{-1}\left( 1+\frac{x^{2}-1}{q-q^{-1}}\,\sigma_i+ 
\frac{x^{2}-1}{1-\nu^{-1}q^{-1} x^{2}}\,\kappa_i\right)\ .\ee
Their classical counterparts (for the Brauer algebras) were found in \cite{Zam}. The elements 
(\ref{bmwbax}) are normalized, $\sigma_i(1)=1$, and satisfy the same unitarity conditions 
(\ref{hunit}). The spectral decomposition of the generator $\sigma_i$ contains three idempotents. 
The basic symmetrizer (the idempotent corresponding to the eigenvalue $q$) is proportional to 
$\sigma_i(q)$. However, $\sigma_i(q^{-1})$ is a mixture of two other idempotents. There are again 
isomorphisms 
$\iota:\, B\!M\!W_{M+1}(q,\nu) \simeq B\!M\!W_{M+1}(-q^{-1},\nu)$,
$\sigma_i\mapsto\sigma_i$. The formula (\ref{bmwbax}) is not invariant under $\iota$.
The basic anti-symmetrizer (the idempotent corresponding to the eigenvalue $-q^{-1}$) is
proportional to $\iota^{-1}(\sigma(x))$ at $x=q$.

\vskip .2cm
Again, the symmetrizers $S_n$, $n=1,\dots,M+1$, are the 
non-zero elements, which satisfy 
\be\lb{desyb} S_1=1\ ,\ S_n\in B\!M\!W_n(q)\subset B\!M\!W_{M+1}(q)\ ,\ee
(\ref{desy'}) and (\ref{desy2}); they exist and are defined uniquely by the conditions 
(\ref{desyb}), (\ref{desy'}) and (\ref{desy2}).

\vskip .2cm
Again, with the knowledge of the baxterized elements, one constructs the symmetrizers 
immediately: they are given by the same formula (\ref{mshu}) and satisfy the same
recurrence (\ref{recm}); the anti-symmetrizers are related to the symmetrizers by the 
isomorphisms $\iota$.

\vskip .2cm
The recurrency (\ref{santis2}) holds for the BMW symmetrizers as well (it was used in 
\cite{WT},\cite{Fior}); it is derived from the baxterized form of the symmetrizers by the
same calculation (\ref{dersan}).

\vskip .2cm
The recurrency relation (\ref{reca}) does not hold for the BMW symmetrizers; the
additive shuffles have to be modified. A version of such modification was suggested
in \cite{HS} and can be derived by a calculation similar to (\ref{santider}). For the 
Hecke algebras the expansions of the products $\sha_{1,n-1}\sha_{1,n-2}\dots\sha_{1,n-k}$ 
contain only reduced words; this is not any more so for the modified shuffles for the BMW 
algebras, the expansions contain similar terms (in a monomial basis, like the one suggested in 
\cite{K}) and the formulas are not as elegant as for the Hecke algebras.

\vskip .2cm
The formula (\ref{syshh}) holds, with the same derivation, for the BMW symmetrizers.

\section{Spectrum of 1-shuffles}\lb{secspec}

Polynomial identities for the multiplicative (for the Hecke and BMW algebras) and additive 
(for the Hecke algebras) 1-shuffles follow, as a by-product, from (\ref{syshh}) and (\ref{sysha}). 
We establish the multiplicities of the eigenvalues in this section. The polynomial identity for 
the additive shuffle was discovered in \cite{Wal} for the symmetric groups and generalized to 
the Hecke algebras in \cite{Lus} with the help of the interpretation of the Hecke algebras in 
terms of flag manifolds over finite fields. The multiplicities of the eigenvalues of the additive 
shuffles were obtained in \cite{DFP} for the symmetric groups. We propose a different approach to 
the calculation of the multiplicities for the Hecke algebras; our method uses the traces of the 
operators of the left multiplication by the additive shuffles.
 
\vskip .2cm
Let $u\in H_n(q)\subset H_m(q)$, $m\geq n$. Denote by $L_u$ the operator of the left 
multiplication by $u$, $L_u:H_m(q)\rightarrow H_m(q)$, 
$L_u(x):=ux$. We denote by ${\mathrm{Tr}}_{_{H_{m}}}(L_u)$ the trace of the operator
$L_u$, considered as the operator on $H_m(q)$.

\paragraph{1.} We start with the multiplicative shuffles. Since
\be\lb{ein} \, ^{^q}\!\Sha_{1,n}S_{n+1}=[n+1]_q!S_{n+1}\ ,\ee 
we find, multiplying (\ref{syshh}) by $(\, ^{^q}\!\Sha_{1,n}-[n+1]_q!)$, the following polynomial 
identity for the multiplicative shuffle
\be (^{^q}\!\Sha_{1,n})^n\,(^{^q}\!\Sha_{1,n}-[n+1]_q!)=0\ ,\ee
which holds for both Hecke and BMW algebras. This is the minimal polynomial, already for the Hecke 
algebras. It is seen without calculations in the Burau representation \cite{Bu} of $H_{n+1}$, 
\be\sigma_j(q^j)\mapsto [j+1]_q\,{\mathrm{Id}}_{j-1}\oplus\left(\begin{matrix}q^{-j}&[j]_q\\
[j]_q&q^j\end{matrix}\right)\oplus [j+1]_q\,{\mathrm{Id}}_{n-j}\ .\ee
If the minimal polynomial is $t^i(t-[n+1]_q!)$ with $i<n$ (the eigenvalue $[n+1]_q!$ is present 
due to (\ref{ein})) then $S_{n+1}$ is proportional to the smaller than $n$ power of 
$\, ^{^q}\!\Sha_{1,n}$.  The matrix of the element $\sigma_j(q^j)$ in the Burau representation
has only one non-zero entry under the main diagonal, on the intersection of {\mbox{$(j+1)$-st}} 
line and $j$-th column. Therefore, the matrix of $^{^q}\!\Sha_{1,n}$ has only one sub-diagonal 
filled with (possibly) non-zeros. However, the matrix of $S_{n+1}$ is proportional to the Hankel 
type matrix $A^i_j:=q^{i+j}$; a smaller than $n$ power of the matrix of $^{^q}\!\Sha_{1,n}$ has 
zero in the very left entry of the bottom line and cannot be equal to the matrix of $S_{n+1}$.

\vskip .2cm
Thus, the element $^{^q}\!\Sha_{1,n}$ is not semi-simple for $n>1$; the semi-simple part of 
$^{^q}\!\Sha_{1,n}$ is $[n+1]_q!S_{n+1}$ and the eigenvalue $[n+1]_q!$ is simple (the rank of 
the projector $L_{S_{n+1}}$ on $H_{n+1}(q)$ is one, because $S_{n+1}\sigma_j=qS_{n+1}$,
$j=1,\dots,n$).

\paragraph{2.} Similarly, multiplying (\ref{sysha}) by $(\sha_{1,n-1}-q^{1-n}[n]_q)$, we find
the following polynomial identity for the additive shuffle
\be\lb{og1}\left(\sha_{1,n-1}- q^{1-n}[n]_q\right) 
\prod_{j=0}^{n-2}\,\left(\sha_{1,n-1}-q^{1-j}[j]_q\right)=0\ . \ee
The $q$-numbers $q^{1-j}[j]_q\equiv 1+q^2+\dots+q^{2j-2}$, $j=1,2,\dots$, are polynomials in 
$q$, linearly independent over ${\mathbb{Z}}$. Therefore, there is a unique integer combination 
$\sum_{j\in\{ 1,2,\dots,n-2,n\} }n_jq^{1-j}[j]_q$, $n_j\in {\mathbb{Z}}$, of these $q$-numbers, 
which is equal to the trace of $L_{\sha_{1,n-1}}$; the coefficients $n_j$ in this combination 
are the multiplicities of the eigenvalues $q^{1-j}[j]_q$, $j>0$. The multiplicity $n_0$ of the 
eigenvalue $0$ is fixed by $\sum n_j={\mathrm{dim}}(H_n(q))\equiv n!$. Thus, the presence 
of the parameter $q$ gives a simple way to calculate the multiplicities. 
 
\vskip .2cm\noindent{\bf Lemma 1.} (i) If $u\in H_j$ then 
\be {\mathrm{Tr}}_{_{H_{j+1}}}(L_u)={\mathrm{Tr}}_{_{H_{j+1}}}(L_{u^{\uparrow 1}})=
(j+1){\mathrm{Tr}}_{_{H_j}}(L_u)\ .\ee
\be {\mathrm{(ii)}}\hspace{2.0cm}{\mathrm{Tr}}_{_{H_{j+1}}}(L_{\sigma_1\sigma_2\dots\sigma_j})
=(q-q^{-1}){\mathrm{Tr}}_{_{H_j}}(L_{\sigma_1\sigma_2\dots\sigma_{j-1}})\ \ ,\ \ j>0
\ .\hspace{.8cm}\ee
\be\lb{dltrsh} {\mathrm{(iii)}}\hspace{2.2cm}
{\mathrm{Tr}}_{_{H_j}}(L_{\sigma_{k-l+1}\dots\sigma_{k-1}\sigma_{k}})
=\frac{j!}{(l+1)!}(q-q^{-1})^l\ \ ,\ \ j>k\geq l\ .\hspace{.7cm}\ee

\vskip .2cm\noindent{\it Proof.} Recall that, as a vector space, 
$H_{j+1}(q)\simeq\oplus_{k=-1}^{j-1}W_k$, where $W_k$ is the vector space 
consisting of elements $v\sigma_j\sigma_{j-1}\dots\sigma_{j-k}$ with $v\in H_j(q)$ (the word 
$\sigma_j\sigma_{j-1}\dots\sigma_{j-k}$ is, by definition, empty for $k=-1$); each 
$W_k$ is canonically isomorphic to $H_j(q)$ as a vector space, the isomorphism is 
$v\sigma_j\sigma_{j-1}\dots\sigma_{j-k}\mapsto v$. The Hecke versions of the automorphism 
$\mathfrak a$ and the anti-automorphism $\mathfrak b$, defined in (\ref{aau}), transform
the above decomposition of $H_{j+1}(q)$ into $H_{j+1}(q)\simeq\oplus_{k=-1}^{j-1}W_k'$, where 
$W_k'$ consists of elements $v\sigma_1\sigma_2\dots\sigma_{k+1}$ with 
$v\in H_j(q)^{\uparrow 1}$ and $H_{j+1}(q)\simeq\oplus_{k=-1}^{j-1}W_k''$, where 
$W_k''$ consists of elements $\sigma_{j-k}\dots\sigma_{j-1}\sigma_jv$ with $v\in H_j(q)$.
 
\vskip .2cm
The operator $L_u$ (respectively, $L_{u^{\uparrow 1}}$) acts in each of the spaces $W_k$
(respectively, $W_k'$) separately and this action commutes with the isomorphisms 
$W_k\simeq H_j(q)$ (respectively, $W_k'\simeq H_j(q)^{\uparrow 1}$). This establishes (i).

\vskip .2cm
In fact, it was enough to find the formula for ${\mathrm{Tr}}_{_{H_{j+1}}}(L_u)$;
${\mathrm{Tr}}_{_{H_{j+1}}}(L_u)={\mathrm{Tr}}_{_{H_{j+1}}}(L_{u^{\uparrow 1}})$
because $u^{\uparrow 1}$ is conjugate to $u$,
$\sigma_1\dots\sigma_ju=u^{\uparrow 1}\sigma_1\dots\sigma_j$, for $u\in H_j(q)$.

\vskip .2cm
(ii) Given a basis $\{ e_\alpha\}$ of a vector space $U$ we say that the vector $e_\alpha$ 
(from the basis) does not contribute to the trace of an operator $X:U\rightarrow U$ if 
$X_\alpha^\alpha=0$ (no summation), where $X_\alpha^\beta$ is the matrix of $X$ in the basis 
$\{ e_\alpha\}$. 

\vskip .2cm
We use the decomposition 
$H_{j+1}(q)\simeq\oplus_{k=-1}^{n-1}W_k''$. The operator
$L_{\sigma_1\sigma_2\dots\sigma_j}$ maps $W_{-1}''$ to $W_{j-1}''$, so vectors from 
$W_{-1}''$ do not contribute to the trace of $L_{\sigma_1\sigma_2\dots\sigma_j}$. For $0\leq k<j$,
\be\lb{cotr}\begin{array}{l} \Bigl(\sigma_1\dots\sigma_j\Bigr)\Bigl(\sigma_{j-k}\dots
\sigma_{j-1}\cdot\sigma_j\Bigr)v=\Bigl(\sigma_{j-k+1}\dots\sigma_j\Bigr)\Bigl(\sigma_1
\dots\sigma_j\Bigr)\sigma_jv\hspace{4cm}\\[1em]
\hspace{1.5cm}=\Bigl(\sigma_{j-k+1}\dots\sigma_j\Bigr)\Bigl(\sigma_1\dots\sigma_{j-1}\Bigr)
\Bigl((q-q^{-1})\sigma_j+1\Bigr)v\\[1em]
\hspace{1.5cm}=(q-q^{-1})\Bigl(\sigma_{j-k+1}\dots\sigma_j\Bigr)\Bigl(\sigma_1\dots\sigma_j
\Bigr)v+\Bigl(\sigma_{j-k+1}\dots\sigma_j\Bigr)\Bigl(\sigma_1\dots\sigma_{j-1}\Bigr)v\\[1em]
\hspace{1.5cm}=(q-q^{-1})\Bigl(\sigma_1\dots\sigma_j\Bigr)\Bigl(\sigma_{j-k}\dots
\sigma_{j-1}\Bigr)v+\Bigl(\sigma_{j-k+1}\dots\sigma_j\Bigr)\Bigl(\sigma_1\dots\sigma_{j-1}
\Bigr)v\ ,\end{array}\ee
the operator $L_{\sigma_1\sigma_2\dots\sigma_j}$ maps $W_k''$ to $W_{j-1}''\oplus W_{k-1}''$.
Therefore, vectors from $W_k''$ do not contribute to the trace of 
$L_{\sigma_1\sigma_2\dots\sigma_j}$ for $k<j-1$. For $k=j-1$, the component 
$L_{\sigma_1\sigma_2\dots\sigma_j}^\diamond$ of the operator $L_{\sigma_1\sigma_2\dots\sigma_j}$, 
which maps $W_{j-1}''$ to $W_{j-1}''$, may have a non-zero trace. This component reads, by 
(\ref{cotr}),
$$L_{\sigma_1\sigma_2\dots\sigma_j}^\diamond(\sigma_1\dots\sigma_jv)=(q-q^{-1})\,
\sigma_1\dots\sigma_jL_{\sigma_1\sigma_2\dots\sigma_{j-1}}(v)\ ,$$ 
and the assertion (ii) follows.

\vskip .2cm
(iii) Follows from (i) and (ii).\hfill$\Box$

\vskip .2cm
By (\ref{dltrsh}), the trace of the operator $L_{\sha_{1,n-1}}$ is
\be {\mathrm{Tr}}_{_{H_n}}(L_{\sha_{1,n-1}})=\sum_{i=0}^{n-1}\frac{n!}{(i+1)!}(q^2-1)^i\ .\ee

{}For the symmetric group ${\mathbb{S}}_n$, the multiplicity of the eigenvalue $j$ of
$L_{\sha_{1,n-1}}$ is the number of permutations in ${\mathbb{S}}_n$ with exactly $j$ fixed 
points \cite{DFP}. Recall that the derangement number $d_n$ (the number of permutations in 
${\mathbb{S}}_n$ without fixed points) is:
\be d_n=n!\,\sum_{i=0}^n\frac{(-1)^i}{i!}\ee
and the number $d_{n,j}$ of permutations in ${\mathbb{S}}_n$ with exactly $j$ fixed points is
\be d_{n,j}=\binc{n}{j} d_{n-j}\equiv\frac{n!}{j!}\sum_{i=0}^{n-j}\frac{(-1)^i}{i!}\ .\ee
{}For generic $q$, the multiplicities of the eigenvalues of $L_{\sha_{1,n-1}}$ are the same as 
for the symmetric group. By construction, $\sum_{j=0}^{n}d_{n,j}=n!$. Thus, to rederive the 
multiplicities we have only to check that $\sum_{j=0}^n d_{n,j}q^{1-j}[j]_q
={\mathrm{Tr}}_{_{H_n}}(L_{\sha_{1,n-1}})$, or, explicitly,
\be\lb{fid}\sum_{j=0}^n\frac{n!}{j!}\sum_{k=0}^{n-j}\frac{(-1)^k}{k!}q^{1-j}[j]_q
=\sum_{i=0}^{n-1}\frac{n!}{(i+1)!}(q^2-1)^i\ .\ee
It is straightforward to verify that both left and right hand sides 
satisfy the same recurrency relation in $n$,
\be f_{n+1}=(n+1)f_n+(q^2-1)^n\ ,\ee
with the same initial condition $f_0=0$, and thereby coincide.

\vskip .5cm
\noindent {\bf Acknowledgements.} The work of A. P. Isaev was partially supported by the 
grants RFBR 08-0100392-a, RFBR-CNRS 07-02-92166-a and RF President Grant N.Sh. 195.2008.2; 
the work of O. V. Ogievetsky was partially supported 
by the ANR project GIMP No.ANR-05-BLAN-0029-01.

\end{document}